\documentclass[twocolumn,dvips]{mtn}
\usepackage{mtnplots}
\usepackage{latexsym}

\begin{document}
\begin{mtnarticle}
\setcounter{equation}{0}
\setcounter{figure}{0}
\setcounter{table}{0}

\title{Numerical Calculations Using Maple: Why \& How?}
\author{
E.V. Corr\^ea Silva\thanks{
     Centro Brasileiro de  Pesquisas F\'\i sicas,
     R. Dr. Xavier Sigaud, 150, Urca, CEP 22290--180,
     Rio de Janeiro, RJ, Brazil. E-mail: ecorrea@cat.cbpf.br},
L.G.S. Duarte\thanks{
     Universidade do Estado do Rio de Janeiro,
     Instituto de F\'{\i}sica, Departamento de  F\'{\i}sica Te\'orica,
     R. S\~ao Francisco Xavier, 524, Maracan\~a, CEP 20550--013,
     Rio de Janeiro, RJ, Brazil. E-mail: lduarte@dft.if.uerj.br}
L.A. da Mota\thanks{
     Idem. E-mail: damota@dft.if.uerj.br}
and J.E.F. Skea\thanks{Idem. E-mail: jimsk@dft.if.uerj.br}
}
\shorttitle{}
\abstract{ The possibility of interaction between Maple and numeric
compiled languages in performing extensive numeric calculations is
exemplified by the {\tt Ndynamics} package, a tool for studying the
(chaotic) behavior of dynamical systems. Programming hints concerning
the construction of {\tt Ndynamics} are presented. The {\tt system}
command, together with the application of the black-box concept, is
used to implement a powerful cooperation between Maple code and some
other numeric language code.}
\keywords{Dynamical systems, chaos, numeric calculations, non-symbolic
compiled languages, black-box}
\maketitle

\section{Introduction}

In this paper we explore the possibility of interaction between Maple
and numeric languages, exemplified by the {\tt Ndynamics} package
\cite{PaperCPC} --- a tool for studying the (chaotic) behavior of
dynamical systems.

Maple interaction facilities may benefit users and programmers of both
symbolic and numeric languages. As far as the implementation of
extensive {\em numerical} calculations is concerned, non-symbolic
compiled languages like {\tt FORTRAN}, {\tt PASCAL} or {\tt C} are
still more popular than any symbolic environment; the existing
collections of efficient programs written in such languages could not
be simply overlooked. 

This paper is organized as follows: first, the concept of
dynamical systems is introduced to the
reader, as well as to the   
concepts of fractal dimension and boundary \cite{ott},
extensively used in the paper; the commands of the package are then briefly
described and an
example of utilization follows. Finally we point out some
useful considerations about programming design, as to the possibility
of interaction between Maple and other languages.

\section{Dynamical Systems}
\label{dynsys}

Systems of ordinary differential equations (or {\em dynamical
systems}) play a central role in a large number of problems in all
areas of scientific research, mainly in physics. Hence the importance
of developing tools to study these systems --- in particular, 
nonlinear systems. A $n$-dimensional dynamical system can be
generically represented by

\begin{equation}
\frac{d\,{\bf X}}{d\,t}={\bf F}({\bf X},t),
\label{syst}
\end{equation}

\noindent where ${\bf X}=\left(X_1(t),X_2(t),\ldots,X_n (t) \right)$
and ${\bf F}=\left(F_1({\bf X},t),\right.$ $ \left. F_2({\bf X},t), \ldots,
F_n ({\bf X},t) \right)$. The variables $X_i(t)$, where $i = 1,2,\ldots,n$,
represent the relevant quantities in some physical model, and $t$ is a
continuous parameter (time, for instance). Each $F_i({\bf X},t)$ is an
arbitrary function of the variables {\bf X} and of the parameter $t$. 

Roughly speaking, ``chaoticity'' means extreme sensitivity to small
changes in the initial conditions. This is the tipical case of
3-dimensional nonlinear dynamical systems, among which chaotic
behavior is a rule, rather than an exception: due to nonlinearity,
small fluctuations in the initial conditions propagate dramatically,
so that two neighbouring inicial conditions may yield orbits with
completely distinct asymptotic behavior.

Chaotic systems have a richer structure than ``well-behaved'' ones,
hence their interest. Strange attractors and repellers, as well as
fractal boundaries are examples of
peculiarities of chaotic systems. A measure of the degree of
chaoticity of a system can be obtained by the evaluation of the
fractal dimension of the boundaries.

In the next section, we are going to elaborate further the concept
of fractals and introduce a method to calculate the associated fractal
dimension.
\section{Fractals}
\label{frac}
\subsection{Fractal Dimension}

The idea of fractal will always lack a precise definition
\cite{mandelbrot}. However, when one refers to a set {\it F} as a
fractal, one typically has the following in mind:

\begin{enumerate}
\item {\it F} has a ``fine structure''; i.e., has complex details on
     arbitrarily small scales;
\item {\it F} is too irregular to be described in the traditional
geometrical language, both locally and globally;
\item often, {\it F} has some form of self-similarity, perhaps
     approximate or statistical;
\item usually, the ``fractal dimension'' of {\it F} (defined in some way)
     is greater than its ``topological dimension'';
\item in most cases of interest, {\it F} is defined in a very simple way,
perhaps recursively.
\end{enumerate}

There are a lot of manners to define the dimension of a fractal. An
important one is called the {\it Hausdorff dimension} \cite{haus}, one
possible generalization of the ``primitive'' notion of dimension.
Suppose we have a hypercube of edge {\it a}, its hypervolume
$V$ being given by
\begin{equation}
\label{dim1}
V = a^d,
\end{equation}

\noindent
where $d$ is the hyperspace dimension.
Let us divide the hypercube into {\it N} hypercubic cells of edge {\it
$\epsilon$}, we  have
\begin{equation}
\label{dim2}
V = N \epsilon^d.
\end{equation}

\noindent
Dividing (\ref{dim2}) by (\ref{dim1}) we have:
\begin{equation}
\label{dim3}
1 = N { \left( {\frac{\epsilon}{a}} \right) }^d
\end{equation}
\noindent
Defining $\delta = {\frac{\epsilon}{a}}$ and expressing
the 
number of
cells as a function of $\delta$ (i.e., $N = N(\delta)$),

\begin{equation}
\label{dim4}
1 = N(\delta) \, \delta^d
\end{equation}

\noindent
Solving (\ref{dim4}) for $d$, we
obtain the standard definition of dimension:
\begin{equation}
\label{dim5}
d = - {\frac{ \ln(N(\delta))}{\ln( \delta )}}
\end{equation}

\noindent Let us now measure the dimension of a Cantor set,
constructed through a recursive procedure pictured in figure
\ref{CantorSet}.

\MaplePlotHeight = 20ex
\MaplePlotWidth= 30ex
\begin{figure} \label{CantorSet}
\mapleplot{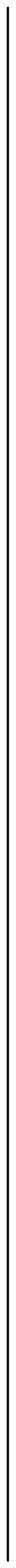}
\mapleplot{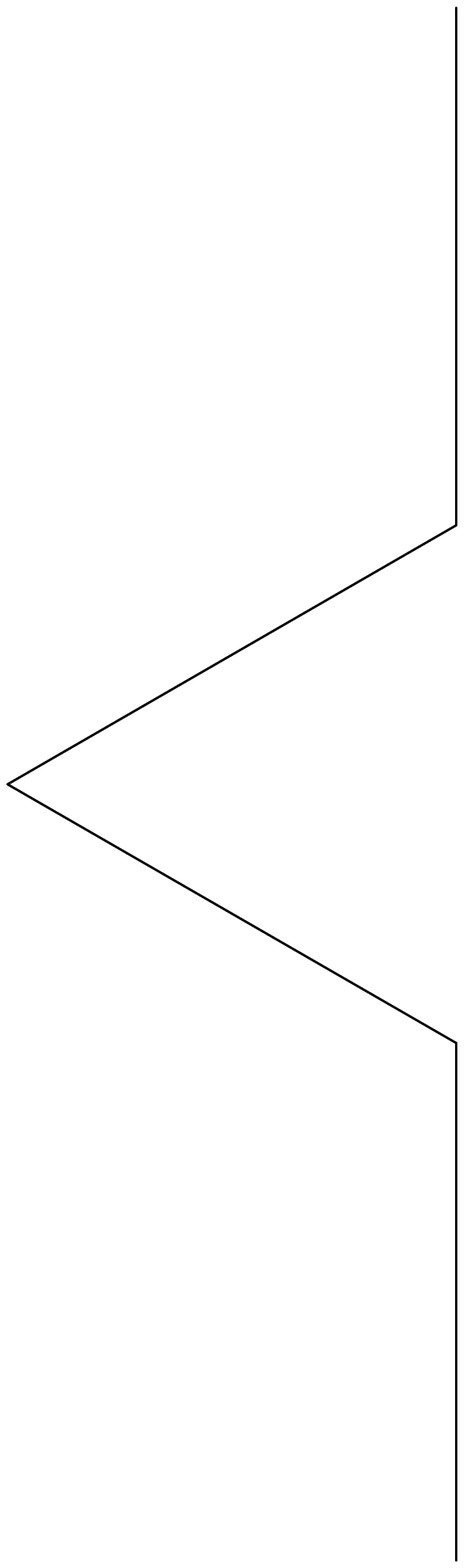}
\mapleplot{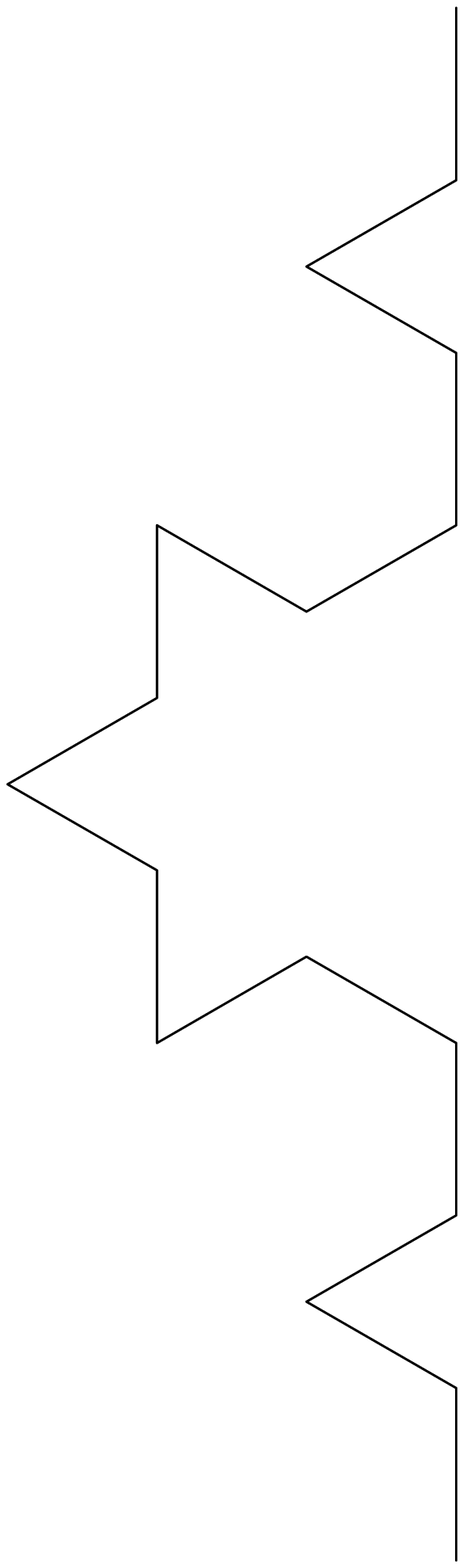}
\mapleplot{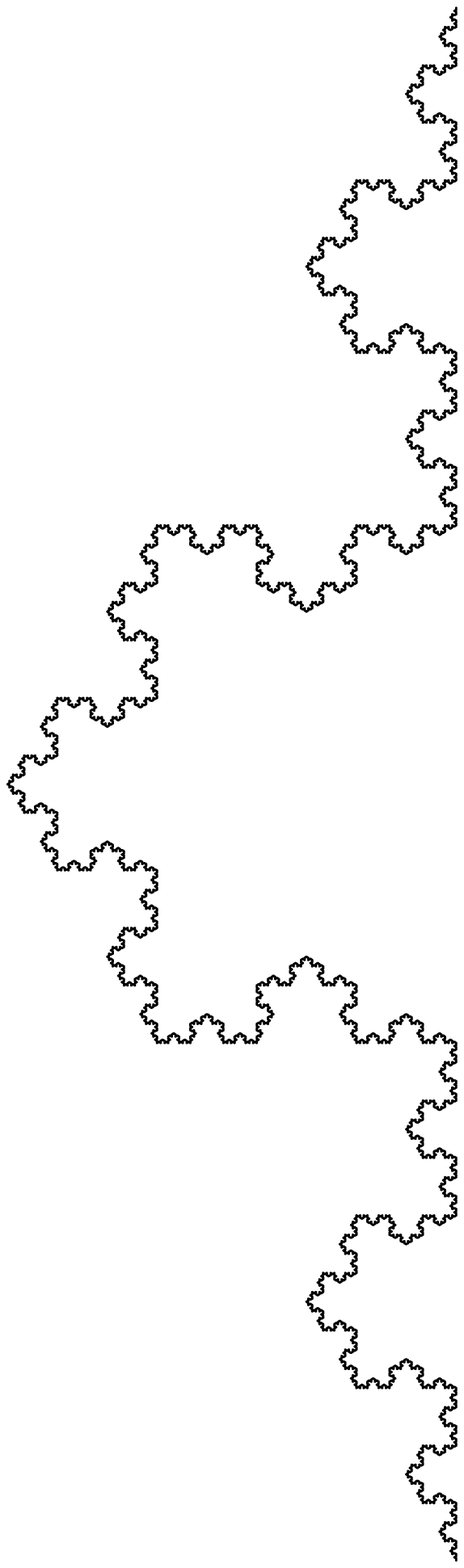}
\caption{Cantor set after 0, 1, 2 and 6 iterations, respectively.}
\end{figure}
          
\MaplePlotHeight = 40 ex 
\MaplePlotWidth = 60 ex  

\noindent In figure 1 we take a segment of length $\epsilon_0$ and divide it into
three parts of equal length $\epsilon_1=\epsilon_0/3$; the middle segment is
then replaced by two segments of length $\epsilon_1$. We
then apply the same procedure to each segment of length $\epsilon_1$,
dividing them into segments of length $\epsilon_2=\epsilon_1/3$, and
so on. After $M$
iterations, we will have a total of $N=4^M$ segments of length
$\epsilon_M = \epsilon_0/3^M$ each. Defining
\begin{equation}
\delta \equiv \frac{\epsilon_M}{\epsilon_0} = \frac{1}{3^M},
\end{equation}
we may write
\begin{equation}
\ln(N(\delta))=-\frac{\ln(\delta)\ \ln(4)}{\ln(3)}.
\end{equation}

\noindent
Taking the limit $M\rightarrow \infty$ (or $\delta \rightarrow
0$), the Hausdorff dimension (\ref{dim5}) of the Cantor set is defined by the limit
\begin{equation}
\label{dim6}
d = -\lim_{\delta \rightarrow 0} {\frac{ \ln(N(\delta))}{\ln \left( \delta \right) }}
\end{equation}

\noindent which, in the case considered, yields
\begin{equation}
\label{dim7}
d = {\frac{ \ln(4)}{\ln \left( {3} \right) }} \sim 1.26
\end{equation}

\subsection{Fractal Dimension of Boundaries}
\label{dimbound}

Let $R$ be a $D$-dimensional finite region --- which we will take to
be a hypercube of edge $a$, without loss of generality --- divided
into $N$ hypercubic ``cells'' of edge $\epsilon$. Let us choose an
arbitrary cell $C$. If there are any pair of points $(P,Q)$ in $C$ so
that the orbits of $P$ and $Q$ have distinct asymptotic behavior, $C$
is said to be a boundary-cell. In the limiting case where $N$ is
infinite ($\epsilon \rightarrow 0$), the union of all boundary-cells 
constitute a {\em boundary}.

Let $N_B$ be the number of boundary-cells in $R$. According to
eq.(\ref{dim6}), the fractal dimension of the boundary \cite{ott} can be inferred
from the behavior of $N_B$ as $\epsilon \rightarrow 0$ or,
equivalently, as $\delta \rightarrow 0$, where $\delta = \epsilon /
a$.

However, the total number of cells
\begin{equation}
N(\delta) = \left( \frac{1}{\delta} \right)^D
\end{equation}
rapidly increases as $\delta \rightarrow 0$, and the computation of
$N_B(\delta)$ results unpractical. An alternative approach is that of
picking up $N^{*}$ random cells in $R$ and counting the number of
boundary-cells $N^{*}_B$. If $N^{*}$ is large enough
to be statiscally meaningful, we expect that
\begin{equation}
\frac{N^{*}_B}{N^{*}} \rightarrow \frac{N_B}{N}.
\end{equation}
Supposing that, indeed,
\begin{equation}
\frac{N^{*}_B}{N^{*}} = \frac{N_B}{N}
\end{equation}
we have
\begin{equation}
\underbrace{- \frac{\ln(N_B)}{\ln \delta}}_{d_B} =
- \frac{\ln({N^{*}_B}/{N^{*}})}{\ln \delta} 
\underbrace{- \frac{\ln(N)}{\ln \delta}.}_{D}
\end{equation}
Defining
\begin{equation}
\alpha \equiv D - d_B = \frac{\ln(N^{*}_B/ N^{*})}{\ln \delta}
\end{equation}
we have
\begin{equation}
\label{alpha}
\ln \left( \frac{N^{*}_B}{N^{*}} \right) = \alpha \ \ln \delta
\end{equation}
Both $\delta$ and $N^{*}_B/N^{*}$ can be measured, allowing us
to determine $\alpha$ and the fractal dimension $d_B$ of
the boundary.
\section{Commands of the Package}
\label{cms}

In this section we present a brief description of the commands of the
{\tt Ndynamics} package. A more complete description can be found in
the on-line help, provided with the 
package \footnote{ http://www.dft.if.uerj.br/symbcomp.html}.

\begin{itemize}
\item {\bf Nsolve} allows the user to specify a system of differential
equations and initial conditions, calculating trajectories in phase
space and generating 2D/3D plots. Also, random initial conditions may
be chosen from a user-defined region of phase space.

\item {\bf View} allows quick and comfortable visual inspection of the
trajectories calculated by {\bf Nsolve}, and the identification of
regions of interest in phase space (e.g., chaos). Zooming in and out
is supported, as well as choosing the variables to be shown (e.g.,
2D or 3D graphs).

\item {\bf Boxcount} performs small random perturbations in the
initial conditions and analyzes their effect upon the evolution of
trajectories. Perturbed initial conditions are chosen from small
regions sized by a user-defined perturbation parameter, around each
umperturbed initial condition. In the neighborhood of a boundary, 
two small but distinct perturbations of the
same initial condition may yield radically different perturbed
trajectories; {\bf Boxcount} determines the fraction of the total
number of initial conditions for which this is so.

\item {\bf Fdimension} manages the execution of {\bf Boxcount} for a
row of values of the perturbation parameter, and analyzes the results 
obtained to evaluate the fractal dimension of the boundary.

\end{itemize}

\section{Examples}
\label{examples}

\subsection{Commands of the package}

{\tt Ndynamics} contains helpful tools for the detection of
interesting regions (such as strange attractors, repellers, and
fractal boundaries) in the phase space of a dynamical system. We will
take the well-known Lorenz system as an example.

The package is loaded by the standard command {\tt with}. The
precision of calculations will be set to 16 digits, using Maple
environment variable {\tt Digits} \footnote{In what follows, the
output of command lines has been ocasionally omitted.}.
\begin{mapleinput}
with(Ndynamics);
Digits := 16;
\end{mapleinput}

\noindent We define the Lorenz system, and assign some values to its
constants:
\begin{mapleinput}
Lorenz := {diff(x(t),t)=sigma*(y(t)-x(t)),
     diff(y(t),t)=-x(t)*z(t)+R*x(t)-y(t),
     diff(z(t),t)=x(t)*y(t)-b*z(t)};
sigma:=10;b:=8/3;R:=28;
\end{mapleinput}

\noindent The command {\tt Nsolve} generates a set of random initial
conditions within a user-defined region in phase space, and then
calculates the orbit of each individual initial condition. The number
of initial conditions to be generated is controled by the global
variable {\tt number\_ic}. In addition, the color of each individual
orbit can be defined by a boolean expression regarding the final condition of
the orbit, to be assigned to the global variable {\tt Colouring}.
According to the value of that expression, an orbit may be plotted in
one out of two predefined colors.

\begin{mapleinput}
number_ic:=8;
Colouring:=x(t)<-4 and x(t) > -11;
\end{mapleinput}

\noindent For the sake of syntax clarity only, the following
assignments are performed: the variable {\tt init\_region} defines
ranges for each dependent variable, thus limiting the region in phase
space for initial conditions; {\tt t\_range} sets the independent
variable range; steps to the independent variable are defined by {\tt
t\_calc\_step} (for calculation purposes) and {\tt t\_plot\_step} (for
plotting purposes).

\begin{mapleinput}
init_region := [x=-0.37717..-0.37716,
     y=0.48685..0.48686, z=-0.29894..-0.29893];
t_range := 0..37;
t_calc_step := .005;
t_plot_step := .01;
\end{mapleinput}

\noindent The command {\tt Nsolve} may now be run. The standard Maple
function {\tt time} is employed to give the reader an estimate
computation time (which is machine-dependent of course).
\begin{mapleinput}
t0:=time():
graph := Nsolve(Lorenz, [init_region,
[t=t_range,t_plot_step]],[x(t),z(t)],
random,method=[rk5C,t_calc_step]):
tf := time()-t0;
\end{mapleinput}
$$ (... some\ ommited\ output\ ...) $$
$$ tf := 123.880$$

\noindent 
We have used the option {\tt method=[rk5C,...]}, i.e., we have used the
`plugged module' in C to perform the numeric calculations.
The result can be plotted, for instance, by using the
command {\tt display} of the standard Maple package {\tt plots}
(previously loaded by {\tt Ndynamics}):
\begin{mapleinput}
display(graph);
\end{mapleinput}
\mapleplot{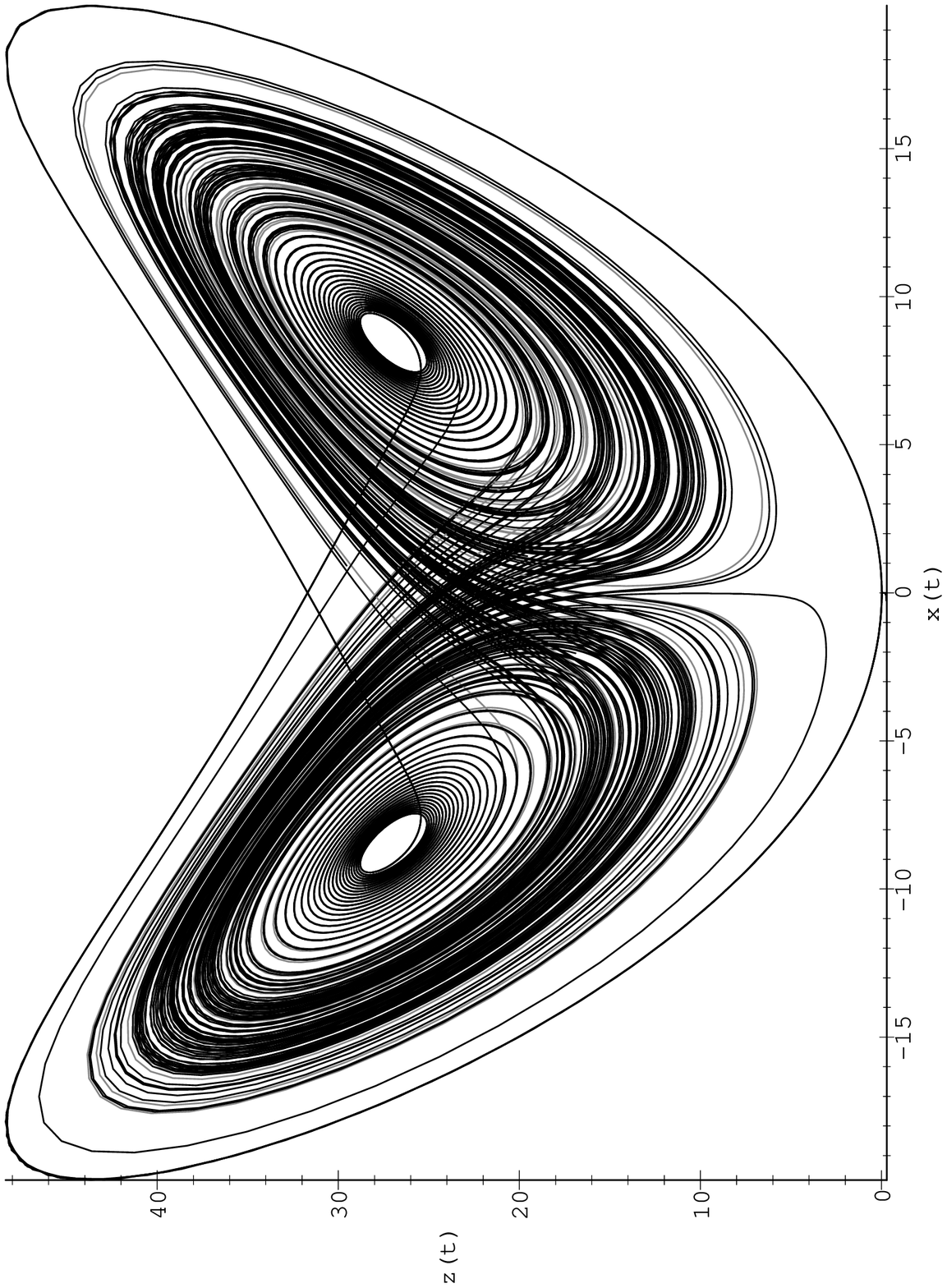} 

\noindent If the user is interested in detailing some particular
region of phase space, the command {\tt View} of {\tt
Ndynamics} is recommended:
\begin{mapleinput}
t0 := time();
View([x=-10..0,z=20..30]); tf := time()-t0;
\end{mapleinput}
\mapleplot{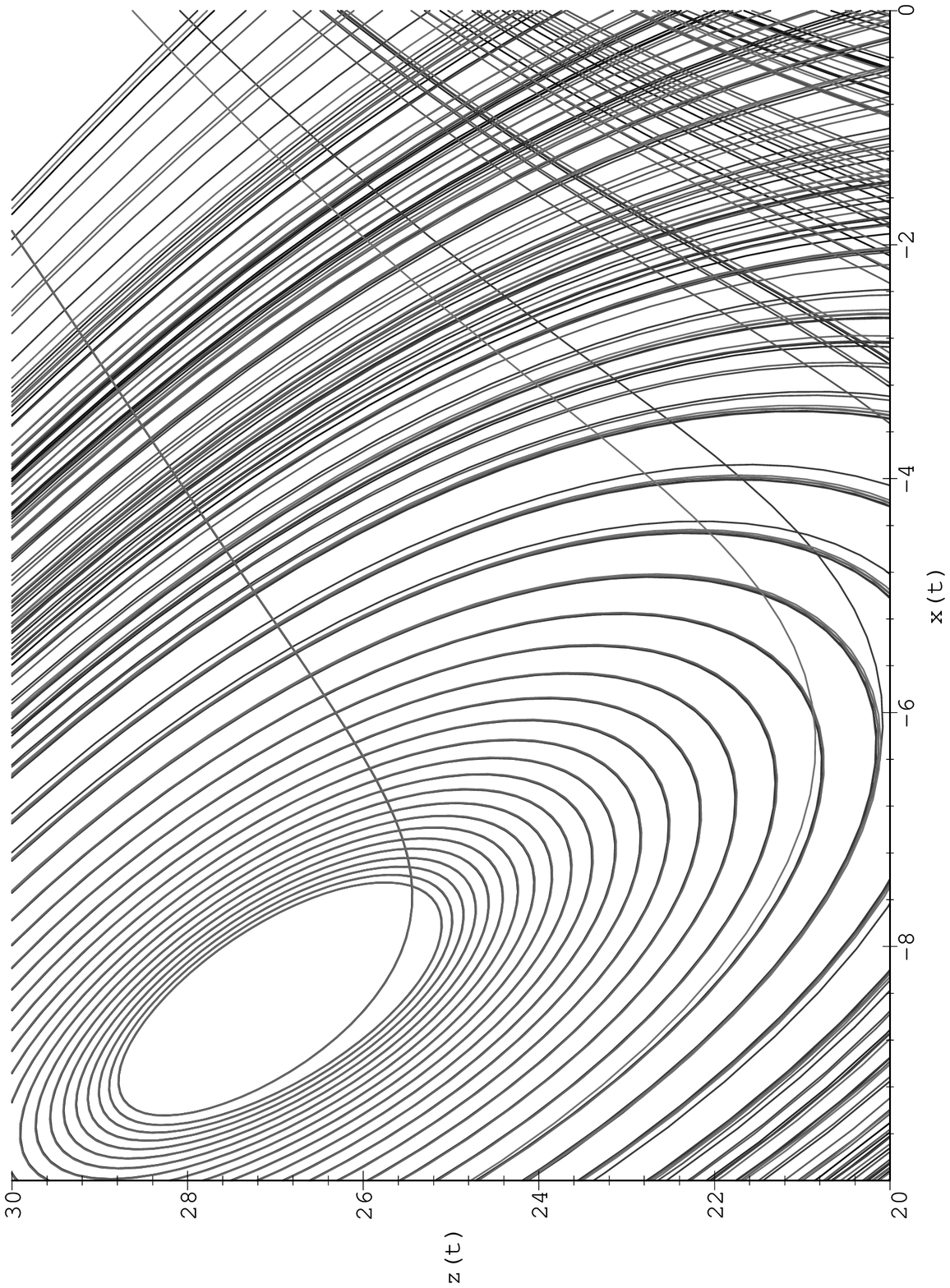} 
$$ tf := 51.569$$

In order to exemplify the use of the commands {\tt Box\-count} and {\tt
Fdimension}, let us slightly modify the Lorenz system parameters:
\begin{mapleinput}
sigma:=10;b:=8/3;R:=20;
\end{mapleinput}
\noindent and choose another region of initial conditions, as well as
other range and steps for the independent variable:
\begin{mapleinput}
init_region2 := [x= -1.001 .. 1.001,
     y = -1.001 .. 1.001,
     z = 21.999 ..22.001];
t_range2 := 0..16;
t_calc_step2 := 0.01;
t_plot_step2 := 0.05;
\end{mapleinput}

\noindent Also, we will choose a much large number of initial
conditions, and a different coloring criteria. Besides, all three
dependent variables $(x,y,z)$ are to be shown in the resulting 3-D
plot:

\begin{mapleinput}
number_ic:=10000;
Colouring:= x(t)<0:
frame := [x(t),y(t),z(t)];
\end{mapleinput}

\noindent Once more {\tt Nsolve} generates random initial
conditions and calculates their orbits:
\begin{mapleinput}
t0 := time():
Nsolve(Lorenz,[init_region2,[t=t_range2,
     t_plot_step2]],frame,initial,random);
tf := time()-t0;
\end{mapleinput}
$$ (...some\ deleted\ output\ ...) $$
$$ tf := 75.725$$

The command {\tt Boxcount} is now able to count the number of
hypercubes over the boundary. Basically, this
command takes perturbed initial conditions (those generated by {\tt
Nsolve}) according to a user-defined parameter, say, {\tt epsilon},
and calculates their orbits. The final value of the independent
variable is represented here by {\tt final\_time}, and the integration
step by {\tt integ\_step}.
(Once more, these variables are assigned just for the sake of clarity
of the {\tt Boxcount} command line.)

\begin{mapleinput}
epsilon := 0.0000002;
final_time := 16;
integ_step := 0.02;
t0 := time():
Boxcount(epsilon,final_time,method=
     [rk5C,integ_step]);
tf := time()-t0;
\end{mapleinput}
$$    reading\ rk5\ output $$
$$ From\ the,\ 10000,\  points\ (that\ were\ testable),\ ,\ 200,$$
$$ of\ them\ were\ close\ to\ the\ boundary. $$
$$ \ [200,\ 10000] $$
$$ tf := 170.560 $$

In its turn, the command {\tt Fdimension} calculates the fractal
dimension of the boundary, by picking up a
given number of distinct values of the perturbation paraneter 
within a user-defined range,
applying {\tt Boxcount} for each case.

\begin{mapleinput}
epsilon_range := 0.0000002..0.000001;
n_epsilons := 5;
t0 := time():
Fdimension(epsilon_range, final_time, n_epsilons,
     method=[rk5C,integ_step]);
tf := time()-t0;
\end{mapleinput}
$$(...some\ deleted\ output ...)$$
$$ Fractal\ dimension\ =\ , 2.213198116174276,\ $$
$$ statistical\ error\ =\, 1.868128109494962,\ \% $$
$$ linear\ correlation\ =\ , .9993933454816381 $$

\mapleplot{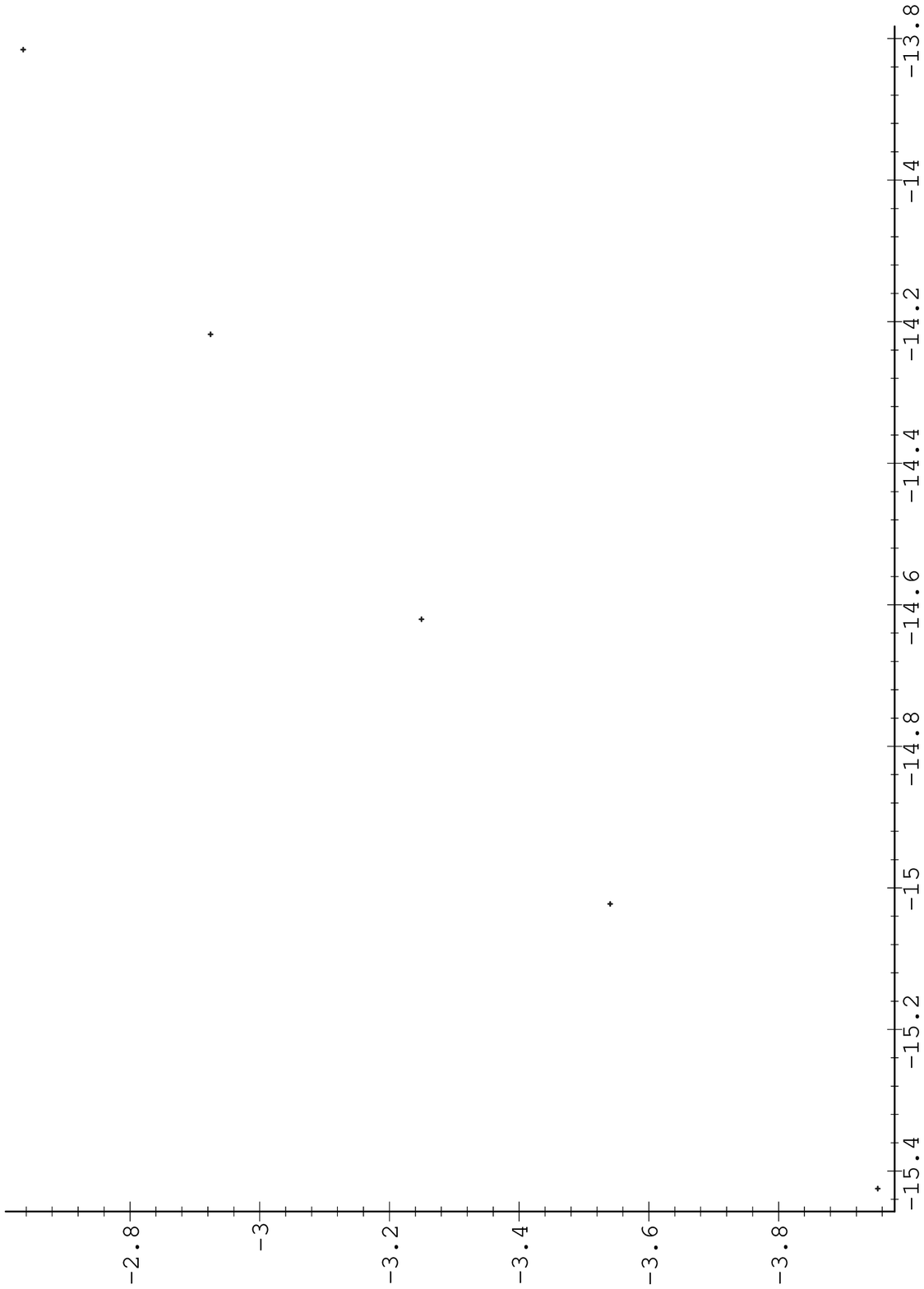} 
$$ tf := 964.940 $$

Where the slope of the straight line above ($\alpha$) is defined by
\ref{alpha}.

\subsection{Performance of the Package}

In the calculations above, the option {\tt method=[rk5C,...]} allowed
us to employ the {\tt Ndynamics} ``plugged module'' --- the file {\tt
rk5.c} --- containing the C source code for the algorithm of the 5-th
order Runge-Kutta method.

To point out the advantages of our hybrid sym\-bol\-ic/nu\-mer\-ic
approach, which uses Maple to manage the source code generation and
compilation for the number crunching, while maintaining Maple's
flexibility, table~\ref{perftab} presents a comparison of the elapsed
time taken to perform the same calculation using the C interface and
performing the entire calculation in Maple, using some of its
numerical integration routines. In order to obtain a fair comparison, particular
care was taken to use methods of the same order and adjusting
parameters to produce solutions with the same precision. The following
procedure was used: for a given set of initial conditions, Maple's
built-in high-precision Taylor Series integrator was used to integrate
the Lorenz system in the interval $0 \leq t \leq 11$, until
convergence was obtained to an accuracy of 13 decimal places;
subsequently Maple's global variable {\tt Digits} was set equal to 13
(higher values of {\tt Digits} seemed to cause problems with Maple's
integrator), and Maple's inbuilt 5th-order Runge-Kutta-Fehlberg
integrator, {\tt rkf45}, was used up to $t=11$. This integration was
found to be accurate to 4 decimal places. The step size in the C
routine was then adjusted to give the same\footnote{in fact, the C
routine was slightly more precise} precision, resulting in a step size
of 0.002. Having matched the precision of both results, the averages of
computation times for both methods was then taken for 200
integrations.

The results below were obtained with Maple~V.4 running in Windows~98
on an AMD-K6 266 with 64Mb of SDRAM. The C compiler used was
Delorie's implementation of GNU's gcc \footnote{http://www.delorie.com},
with level~3 optimization.

\begin{table}[htb]
\begin{center}
\begin{tabular}{|c|c|c|}
\hline 
{\bf Integration}  & {\bf C-interface} & {\bf Maple} \\
{\bf Method}  & {\bf rk5C} &{\bf rkf45} \\
\hline 
Seconds per Trajectory  & 0.1 & 7.0 \\
\hline
Ratio to the Fastest Case  & 1 & 70 \\
\hline 
\end{tabular}
\caption{Comparative performances of Maple's inbuilt numerical integrators
and {\tt Ndynamics} C interface, for the Lorenz system.}
\label{perftab}
\end{center}
\end{table}

\section{Programming and Design \\ Considerations}
\label{progdesing}

Building computer programs demands methodical approach, if time and
effort of construction are to be saved, and coherence to be kept.
Furthermore, debugging and modification of a program relates directly
to how apparent the {\em purpose} of each piece of code turns out.
{\em Software design} is, in itself, a rich subject and field of
research, barely touched here. Yet, we hope the reader may profit from
reading our present considerations.

\subsection{Interface with Compiled Languages}
\label{inter}

The numerical integration of a system of differential equations
exemplifies a tipical situation involving extensive numerical
calculations. As far as the efficiency (basically, memory usage and
processing time) of numerical algorithms is concerned, general-purpose
symbolic computing systems (based on {\em interpreted} languages)
generally offer less attractive means of implementation than
non-symbolic {\em compiled} languages such as C, FORTRAN, and PASCAL.
The long history of non-symbolic languages has already produced a
large collection of efficient programs for a wide variety of specific
problems --- well-known examples are the NAG library of FORTRAN
programs, and the collection of C programs in \cite{NumericalRecipes}.
Using a symbolic environment for numerical calculations, discarding
existing solutions (based on  compiled languages) might sound like
an attempt to ``reinvent the wheel''.

The point is that symbolic languages provide flexibility and
interactivity which lacks in ``numeric-oriented'' languages. In other
words, symbolic code and data structures tend to be more friendly and versatile than
numeric ones. For instance, in our specific
case, the identification of chaotic regions in phase space is made
much comfortable by the interactive user-interface of Maple. (Changing
parameters and checking out results is faster and easier.)

In one aspect, the question is analogue to ``Why using a high-level
language instead of Assembly ?'' No matter how fast Assembly programs
may be, the kind of details involved in their construction (most of
them totally unrelated to the problem at hand) certainly discourage
its use. The more a given language allows its programmer to focus on
``problem-related'' details and {\em ignore} the rest, the most
suitable the language becomes to find a solution to the problem. The
very nature of ever-evolving problems in scientifical research
recommends the use of ``higher-level'' languages.

It would be highly desirable to combine the features of symbolic and
numeric languages, for the task in view; that was one of the
guidelines to the construction of our package. ``High-level'' tasks,
so to speak, are handled by code written in Maple, whereas numerical
calculations may be (optionally) performed by a piece of code written
in some other language. The present implementation of the program
contains an extension written in C; future versions will allow the
user to simply ``plug'' his own (numerical) code. The implementation
of these ideas in Maple was achieved through the utilization of the
standard command {\tt system} \footnote{Please refer to the Maple
on-line help, by typing {\tt ?system}.} --- which allows one to send
commands directly to the operational system --- as well as the
possibility of exchanging data through files.

In the current implementation of {\tt Ndynamics}, the command {\tt
system} is used for:
\begin{itemize}
\item{ compiling \footnote{Currently, the DJGPP-32 GNU compiler is
used. Please refer to {\tt Ndynamics} on-line help for more details.}
two source code files: {\tt derivs.c}, containing the definition of the
dynamical system; and {\tt rk5.c}, containing the code for the 5-th
order Runge-Kutta method. The executable file {\tt rk5.exe} is then
produced.}

\item{executing {\tt rk5.exe} to perform calculations. The files {\tt
rk5.in} and {\tt rk5.out} are, respectively, the input and output
files for {\tt rk5.exe}.}

\end{itemize}

\noindent All these operations are performed automatically, entirely
from within Maple. The user needs not to be aware of any such
processes.

The attentive Maple programmer will notice that the {\it utilization}
of {\tt rk5.exe} --- assuming that {\tt rk5.in} and {\tt rk5.out} have
the proper formats --- does not depend whether its original source
code was in C or any other language. This is an example of a ``
black-box'', a concept to be explored in what follows.

\subsection{Black-Boxes and Parameter Passing Techniques}
\label{programing}

A well-know concept in software design is that of a black-box
\cite{YourdonOld}. For our present purposes, it can be
pictured as a Maple procedure or an executable file which can be {\em
used} without any assumption on its internal logic: all we need to
know relates either to input and output data (meaning, form, values)
or to {\em what} the procedure does (rather than {\em how} it is
done). A good example of a black-box is the standard Maple function
{\tt sin}: all one needs to know to use {\tt sin} is that it accepts a
single input argument (a valid algebraic expression, assumed to be
expressed in radians) and that its output is, also, an algebraic
expression. Using the function {\tt sin} does not require the
knowledge of details concerning actual computations.

Ideal packages, made out of black-box-like procedures, would be most
easily understood, and hence most easily tested, debugged and
modified.

Exchange of data among procedures is a crucial aspect in the
construction of a package. A Maple procedure may exchange data by
\begin{itemize}
\item using global or environment variables (input and output);
\item using the result of its last executed line (output);
\item using the RETURN command (output);
\item reading from or writing to a file (input and output);
\item using its parameters or arguments (input and output).
\end{itemize}
\noindent

\noindent Of all these ways, the last one appears to us as the most
interesting, as far the implementation of black-boxes
is concerned. Using arguments of a procedure as a means of
conveying {\em input} data is a common practice,
needing no further comments.
We would like to invite the reader to attend more closely to the
possibility of passing {\em output} data through the arguments of a
procedure --- which is also a known technique
\cite{Heck,FirstLeaves} found, for instance, in the standard
commands {\tt iquo} and {\tt member}:
\begin{mapleinput}
iquo(11,2,'r');
 # r is an output parameter, to be
 # internally assigned to the
 # remainder of the integer division
 # 11/2
\end{mapleinput}
\begin{maplelatex}
\[ 5 \]
\end{maplelatex}
\begin{mapleinput}
 r;
\end{mapleinput}
\begin{maplelatex}
\[ 1 \]
\end{maplelatex}
\begin{mapleinput}
member(x,[y,y,x,y],'pos');
   # pos is an output parameter, to be
   # assigned to the position of the
   # first occurence of x in [y,y,x,y]
\end{mapleinput}
$$ true $$
\begin{mapleinput}
pos;
\end{mapleinput}
$$ 3 $$

\noindent Another example is the command {\tt irem}, similar to {\tt
iquo}. Let us make our point clear with the help of a ``toy
procedure'' {\tt f}:
\begin{mapleinput}
f := proc(a::posint,sm_,pr_)
local i;
sm_ := convert([seq(i,i=1..a)],`+`);
pr_ := convert([seq(i,i=1..a)],`*`);
NULL;
end:
\end{mapleinput}

\noindent Given a positive integer {\tt a}, we compute the sum
($1 +2+... +\ ${\tt a}) and the product ($1 * 2 * ... \ * $ {\tt a}),
which are assigned to the parameters {\tt sm\_} and {\tt pr\_}\ ,
respectively. An example of usage is
\begin{mapleinput}
f(5,'s,p'); # null output
s,p;
\end{mapleinput}
\begin{maplelatex}
\[ 15,\ 120 \]
\end{maplelatex}

\noindent Now compare the procedure {\tt f} to a slightly (but
significantly) modified version {\tt g}, where output parameters are
not used:
\begin{mapleinput}
g := proc(a::posint)
local i, s, p;
s := convert([seq(i,i=1..a)],`+`);
p := convert([seq(i,i=1..a)],`*`);
s,p;
end:
\end{mapleinput}

\noindent This version employs the ``last executed line'' technique.
It would run as
\begin{mapleinput}
g(5);
\end{mapleinput}
\begin{maplelatex}
\[ 15,\ 120 \]
\end{maplelatex}

\noindent What is the practical difference between the two procedures?
The reader should notice that, in order to interpret (and evaluate the
correctness of) the output of {\tt g}, he must hold the {\it extra}
information that the {\em first} operand of the output sequence
corresponds to $1+2+3+4+5$, and the {\em second} one corresponds to
$1*2*3*4*5$ --- i.e., the interpretation of the output depends on its
particular (and arbitrary) form. In this example, it is a sequence,
but it might as well be a list, or even a table. Were it a sequence or
a list, the first operand might be the product, instead of the sum ...
Well, whatever the choice as to the output form--- which {\it has to}
be made ---, it would certainly have little to do with the actual
purpose of the procedure (i.e., calculating sums and products), and
thus would be superfluous detail. A procedure calling {\tt g} (and its
poor programmer) would have to know all about these details, though;
it (or he, or she) would be free from such concern, if {\tt f} were
used instead. The procedure {\tt f} can be used as a black-box,
whereas {\tt g} can't. 

When it comes to developing a large number of integrated Maple
procedures, the black-box is a valuable concept to keep in mind, and
so is the ``output parameter'' technique. That doesn't mean that all
procedures should be written in that way, though. The ``last line''
and ``RETURN'' techniques, for instance, seem to be the most suitable
to user-level procedures. The ``output parameter'' technique, in its
turn, fits best on internal procedures, where all the ``hard work'' is
carried out.

\section{Conclusion}
\label{conclusion}

The powerful combination of Maple flexibility and interactivity with
the high speed of numeric languages has been used in the
implementation of {\tt Ndynamics}, a package for the numerical study
of (chaotic) behaviour in dynamical systems. A more detailed
description of the functions and physical applications of {\tt
Ndynamics} can be found on \cite{PaperCPC}.

We have also presented some of the programming hints and strategies
used in its construction, hoping to have shown that there can be a a
stronger link between users (and programmers) of two kinds: those who
are interested in heavy numerical calculations, and those whose
interest lies in symbolic computations. 

Part of the current implementation of {\tt Ndynamics} has been written
in C and has been ``plugged'', so to speak, into the main package.
Currently, the user himself can write and use his own integration
routine in C, as long as the same form of input and output data is
kept \footnote{For more detailed instructions, please refer to the
on-line documentation.}. Future versions of {\tt Ndynamics} will allow
the user to ``plug and play'' his own numerical piece of code and run
the compiler corresponding to the language of his preference, provided
that the input and output data follow the proper conventions.

\section{Biographies}

\begin{itemize}

\item {\bf E.V. Corr\^ea Silva} is currently a Ph.D. student of
Physics (Quantum Field Theories) at the Centro Brasileiro de Pesquisas
F\'{\i}sicas, in Rio de Janeiro, Brazil. He obtained his Master degree
in Physics at the Universidade Federal do Rio de Janeiro, in 1994. His
research interests include the Quantum Hall Effect, Quantum Field
Theory and Chaos, as well as System Design, Artificial Intelligence,
Computational Physics and its relation to the teaching of Physics.

\item {\bf L.G.S.Duarte} is currently a member of staff at the Universidade Estadual
do Rio de Janeiro (UERJ). He obtained his PhD in Physics at the Universidade
Federal do Rio de Janeiro, in 1997. His interests concentrate on mathematical
physics and he has recently participated on a project that lead to the
creation of a differential equation solver, writen on Maple, that was finally
introduced on the comercial release number 5 of such mathematical package.

\item {\bf L.A.C.P.da Mota} is currently a member of staff at the Universidade
Estadual do Rio de Janeiro (UERJ). He obtained his PhD in Physics at the
University of Oxford, U.K, in 1993. He then held an one year position as a
visiting researcher at the same University. His interests were then Elementary
Particle Phenomenology. Most recently, his work has concentrated on
mathematical physics and participated on the same project (mentioned on the
biography above).

\item {\bf J.E.F.Skea} graduated from the University of Glasgow with a joint honours
degree in Physics and Astronomy. He first worked with computer algebra as
a research tool at the University of Sussex where, under the supervision
of Roger Tayler and John Barrow, he completed his Ph.D. in the areas of
general relativity and cosmology. Subsequently he worked for 3 years as
a postdoc at the School of Mathematical Sciences of Queen Mary College
London with Malcolm MacCallum, mainly on the CA systems SHEEP and REDUCE.
He left science for a couple of years to work in Science Policy with Ben
Martin at SPRU in the University of Sussex, on the analysis of performance
indicators, but returned to physics for his third postdoc, this time at
CBPF (Brazilian Centre for Physics Research) in Rio.  He moved back to
Scotland when he finally got a teaching position in the Department of
Mathematics of the University of Aberdeen, where he initiated a course on
the theory of computer algebra, using Maple in practical sessions. In
1996 he moved to the Theoretical Physics department of UERJ (Rio State
University) where he seems to have settled down ... for now.

\end{itemize}

\end{mtnarticle} 
\end{document}